\documentclass[12pt]{article}
\usepackage{eurosym}

\usepackage[latin1]{inputenc}
\usepackage{amsmath}
\usepackage{latexsym}
\usepackage{graphicx}

\usepackage{amsfonts,amssymb}
\usepackage{tikz}
\usetikzlibrary{automata,arrows,positioning,calc}

\usepackage{epsfig}
\usepackage{supertabular}

\newcommand\be{\begin{equation}}
\newcommand\ee{\end{equation}}
\newcommand\bes{\begin{equation*}}
\newcommand\ees{\end{equation*}}

\def\n0{\mathbb{N}_0}

\def\A{\mathcal{A}}

\def\1{\mathbf{1}}
\def\Elim{E_{\infty}}

\def\cal{\mathcal}

\title{On certain families of planar patterns and fractals}

\author{Ligia L. Cristea \thanks{The author is supported by the Austrian Science Fund (FWF), 
Project P27050-N26, and by the Austrian Science Fund (FWF) Project F5508-N26, which is part of the Special Research Program ``Quasi-Monte Carlo Methods: Theory and Applications''. } 
\\Karl-Franzens-Universit\"at Graz\\ Institut f\"ur Mathematik und Wissenschaftliches Rechnen
\\Heinrichstrasse 36, 8010 Graz,Austria\\ \emph{strublistea@gmail.com} }

\begin{document}

\maketitle

\emph{This paper is dedicated to Professor Robert F. Tichy on the occasion of his 60$^{th}$ anniversary.}
\\[0.1cm]


\abstract{This survey article is dedicated to some families of fractals that were introduced and studied during the last decade, more precisely, families of Sierpi\'nski carpets: limit net sets, generalised Sierpi\'nski carpets and labyrinth fractals. We give a unifying approach of these fractals and several of their topological and geometrical properties, by using the framework of planar patterns.}

\section{Introduction}\label{sec:introduction}

{\em Sierpi\'nski carpets} are self-similar fractals in the plane that originate from the classical Sierpi\'nski carpet \cite{Mandelbrot1983, Sierpinski}. Sierpi\'nski carpets are constructed by dividing the unit square into $m \times m$ congruent smaller subsquares of which $m_0$ squares are cut out together with their boundary, and then taking the closure. The resulting pattern is called the \emph{generator} of the 
Sierpi\'nski carpet. At each step of the iterative construction
this procedure is applied to all remaining squares, and, repeating this construction ad infinitum, the resulting object is a
fractal of Hausdorff dimension $\frac{\log(m^2 -m_0)}{\log(m)}$, called a
\emph{Sierpi\'nski carpet} \cite{the_pore_structure_2001}. Figure 
\ref{fig:laby patterns krabbe} shows the first two steps of the interative construction of a Sierpi\'nski carpet.

These fractals can also be defined as attractors of IFS (for Iterated Functions Systems we refer, e.g.,  to the books of Falconer \cite{Falconerbook1990, Falconerbook1997} and  Barnsley \cite{barnsley,barnsley_superfractals}), and occur in several branches of mathematics. In particular, their geometric and topological properties gained a lot of interest, see, e.g., Whyburn~\cite{MR0099638}, Curtis and Fort~\cite{MR0105103}, McMullen~\cite{McMullen_general_carpets_1984}, Bandt and Mubarak \cite{BandtMubarak_2004}, Lau et al.~\cite{fractal_squares_hongkong}. During the last decades Sierpi\'nski carpets have been used, e.g., as models for
porous materials \cite{the_pore_structure_2001, tarafdar_modelling_porous_structures_2001}.

\emph{Limit net sets, generalised Sierpi\'nski carpets  and labyrinth fractals} are families of Sierpi\'nski carpets that were introduced and studied by Cristea and Steinsky \cite{netsets,connected_general_carpets,totally_disco,laby4x4,generallaby} and some of the results were extended in recent research \cite{mixlaby, notearcsmixlaby}  to even more general fractal objects called \emph{mixed labyrinth fractals}. Studying these objects is  of interest not just for mathematics, but also for research in physics, where some of the results have already been used, e.g., \cite{laby_fizicieni1, laby_fizicieni2, PotapovZhang_oct2016}.

In this paper we present  results on topological and geometrical properties that were obtained for the three families of Sierpi\'nski carpets mentioned above, such as connectedness or lengths of arcs in these fractals, everything being done under a combinatorial frame, where the combinatorial character of the problems comes from the combinatorics of the generator(s) of the carpet: the pattern(s).

Although originally \emph{net sets} and \emph{limit net sets} were defined and constructed by means of net matrices \cite{netsets}, and the \emph{labyrinth fractals} by using \emph{labyrinth sets} \cite{laby4x4,generallaby}, throughout this paper we give a unifying approach of all the families of carpets mentioned above by means of \emph{patterns}, as it was done in the case of the \emph{generalised Sierpi\'nski carpets} that were studied \cite{connected_general_carpets, totally_disco} after the other mentioned carpets, and in more recent work, for mixed labyrinth fractals \cite{mixlaby,notearcsmixlaby}. 

 Graph directed constructions, see, e.g., \cite{MauldinWilliams_graphdirected_1988}, GDMS (Graph Directed Markov Systems, see, e.g., 
\cite{MauldinUrbanski_bookGDMS_2003}), and random fractals \cite{Falcrandomfr1986, MauldinWilliams_random_1986} also offer framewors for studying the objects that occur along this paper.
Finally, we mention that there are recent results and ongoing research on $V$-variable fractals, see e.g. \cite{freiberghamblyhutchinson},
 and several of the fractals studied and mentioned in this section can be approached within the frame of $V$-variable fractals.
For $V$-variable fractals and superfractals we also refer to Barnsley's book \cite{barnsley_superfractals}.

Let us now give a short outline of the paper. In Section \ref{sec:patterns} we define planar patterns and the graph associated to a planar pattern. Section \ref{sec:limit_net_sets} is dedicated to net sets and limit net sets. In Section \ref{sec:GSC} we briefly present recent research on generalised Sierpi\'nski carpets. Section \ref{sec:laby} deals with results about self-similar and mixed labyrinth fractals and also refers to very recent results. Finally, Section \ref{sec:conclusions} is dedicated to conclusions and final remarks of the survey.

\section{Planar patterns and Sierpi\'nski carpets}
\label{sec:patterns}


First, let us recall the definition of a pattern, as it is given in some of the above mentioned papers \cite{connected_general_carpets, totally_disco}.
Let $x,y,q\in [0,1]$ such that $Q=[x,x+q]\times [y,y+q]\subseteq [0,1]\times [0,1]$. Then for any point $(z_x,z_y)\in[0,1]\times [0,1]$ we define the function $\displaystyle P_Q(z_x,z_y)=(q z_x+x,q z_y+y)$.

Let $m\ge 1$. For the integers $i,j$ with $0\le i,j\le m-1,$ let $S_{i,j}^{m}=\{(x,y)\mid \frac{i}{m}\le x \le \frac{i+1}{m} \mbox{ and } \frac{j}{m}\le y \le \frac{j+1}{m} \}$, and  
${\cal S}_m=\{S_{i,j}^{m}\mid 0\le i\le m-1 \mbox{ and } 0\le j\le m-1 \}$. 
We call any nonempty $\A \subseteq {\cal S}_m$ an  $m \times m$ \emph{pattern} or, in short, $m$-\emph{pattern}.

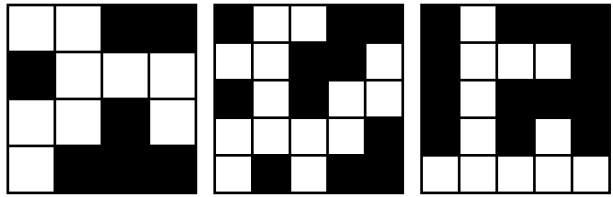
\begin{figure}[hhhh]\label{A1A2}
\begin{center}
\begin{tikzpicture}[scale=.25]
\draw[line width=1pt] (0,0) rectangle (10,10);
\draw[line width=1pt] (2.5, 0) -- (2.5,10);
\draw[line width=1pt] (5, 0) -- (5,10);
\draw[line width=1pt] (7.5, 0) -- (7.5,10);
\draw[line width=1pt] (0, 2.5) -- (10,2.5);
\draw[line width=1pt] (0, 5) -- (10,5);
\draw[line width=1pt] (0, 7.5) -- (10,7.5);
\filldraw[fill=black, draw=black] (2.5,0) rectangle (5, 2.5);
\filldraw[fill=black, draw=black] (5,0) rectangle (7.5, 2.5);
\filldraw[fill=black, draw=black] (7.5,0) rectangle (10, 2.5);
\filldraw[fill=black, draw=black] (5,2.5) rectangle (7.5, 5);
\filldraw[fill=black, draw=black] (0,5) rectangle (2.5, 7.5);
\filldraw[fill=black, draw=black] (5,7.5) rectangle (7.5, 10);
\filldraw[fill=black, draw=black] (7.5,7.5) rectangle (10, 10);
\draw[line width=1pt] (11,0) rectangle (21,10);
\draw[line width=1pt] (13, 0) -- (13,10);
\draw[line width=1pt] (15, 0) -- (15,10);
\draw[line width=1pt] (17, 0) -- (17,10);
\draw[line width=1pt] (19, 0) -- (19,10);
\draw[line width=1pt] (11, 2) -- (21,2);
\draw[line width=1pt] (11, 4) -- (21,4);
\draw[line width=1pt] (11, 6) -- (21,6);
\draw[line width=1pt] (11, 8) -- (21,8);
\filldraw[fill=black, draw=black] (13,0) rectangle (15, 2);
\filldraw[fill=black, draw=black] (17,0) rectangle (19, 2);
\filldraw[fill=black, draw=black] (19,0) rectangle (21, 2);
\filldraw[fill=black, draw=black] (19,2) rectangle (21, 4);
\filldraw[fill=black, draw=black] (11,4) rectangle (13, 6);
\filldraw[fill=black, draw=black] (15,4) rectangle (17, 6);
\filldraw[fill=black, draw=black] (15,6) rectangle (17, 8);
\filldraw[fill=black, draw=black] (17,6) rectangle (19, 8);
\filldraw[fill=black, draw=black] (11,8) rectangle (13, 10);
\filldraw[fill=black, draw=black] (17,8) rectangle (19, 10 );
\filldraw[fill=black, draw=black] (19,8) rectangle (21, 10);
\draw[line width=1pt] (22,0) rectangle (32,10);
\draw[line width=1pt] (24, 0) -- (24,10);
\draw[line width=1pt] (26, 0) -- (26,10);
\draw[line width=1pt] (28, 0) -- (28,10);
\draw[line width=1pt] (30, 0) -- (30,10);
\draw[line width=1pt] (22, 2) -- (32,2);
\draw[line width=1pt] (22, 4) -- (32,4);
\draw[line width=1pt] (22, 6) -- (32,6);
\draw[line width=1pt] (22, 8) -- (32,8);
\filldraw[fill=black, draw=black] (22,2) rectangle (24, 10);
\filldraw[fill=black, draw=black] (26,2) rectangle (28, 6);
\filldraw[fill=black, draw=black] (28,4) rectangle (30, 6);
\filldraw[fill=black, draw=black] (30,2) rectangle (32, 8);
\filldraw[fill=black, draw=black] (26,8) rectangle (32, 10);

\end{tikzpicture}
\caption{Three patterns, ${\cal A}_1$ (a $4$-pattern), ${\cal A}_2$ (a $5$-pattern)  and ${\cal A}_3$ (a $5$-pattern) }
\label{fig:A1A2A3}
\end{center}
\end{figure}

In Figure \ref{fig:A1A2A3} we show three such patterns that are in particular also labyrinth patterns, which we define in Section \ref{sec:laby}. We mention that throughout this paper we think of the black regions in the figures as being ``cut out'' at the corresponding step, and subsequently the closure (with respect to the topology induced by the Euclidean metric in the plane) of the remainder set is taken.
                                                                                
All families of fractals that we present in this paper can be constructed by means of patterns, and in each case we use a iterative construction, analogous to that described in Section \ref{sec:introduction} for a Sierpi\'nski carpet.

For any pattern ${\cal A}\subseteq {{\cal S}_m}$, we define the graph $\mathcal{G}({\cal A})\equiv (\mathcal{V}(\mathcal{G}({\cal A})),\mathcal{E}(\mathcal{G}({\cal A})))$ to be \emph{the graph of ${\cal A}$}, i.e., the graph whose vertices are the (closed) white squares in ${\cal A}$, i.e.,  $\mathcal{V}(\mathcal{G}({\cal A}))={\cal A}$ and whose set of edges $\mathcal{E}(\mathcal{G}({\cal A}))$ consists of the unordered pairs of white squares, that share a common side.

\section{Limit net sets}\label{sec:limit_net_sets}

\par \emph{Net set}  and \emph{limit net set} are new concepts developed in 
\cite{netsets}, based, on the one hand, on the observation that various porous materials present holes that
at each scale are ``evenly'' distributed,
 and, on the other hand, on the distribution properties of $(t,m,s)$-nets, that are well distributed point sets in the unit cube, for more details see \cite{Nied.netsMhM}.

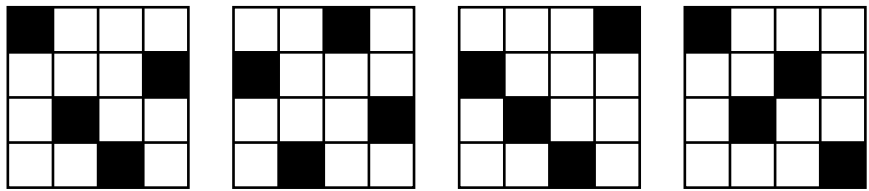
\begin{figure}[hhhh]
\begin{center}
\begin{tikzpicture}[scale=.30]
\draw[line width=1pt] (0,0) rectangle (8,8);
\draw[line width=1pt] (2, 0) -- (2,8);
\draw[line width=1pt] (4, 0) -- (4,8);
\draw[line width=1pt] (6, 0) -- (6,8);
\draw[line width=1pt] (0, 2) -- (8,2);
\draw[line width=1pt] (0, 4) -- (8,4);
\draw[line width=1pt] (0, 6) -- (8,6);
\filldraw[fill=black, draw=black] (0,6) rectangle (2, 8);
\filldraw[fill=black, draw=black] (2,2) rectangle (4, 4);
\filldraw[fill=black, draw=black] (4,0) rectangle (6, 2);
\filldraw[fill=black, draw=black] (6,4) rectangle (8, 6);
\draw[line width=1pt] (10,0) rectangle (18,8);
\draw[line width=1pt] (12, 0) -- (12,8);
\draw[line width=1pt] (14, 0) -- (14,8);
\draw[line width=1pt] (16, 0) -- (16,8);
\draw[line width=1pt] (10, 2) -- (18,2);
\draw[line width=1pt] (10, 4) -- (18,4);
\draw[line width=1pt] (10, 6) -- (18,6);
\filldraw[fill=black, draw=black] (10,4) rectangle (12, 6);
\filldraw[fill=black, draw=black] (12,0) rectangle (14, 2);
\filldraw[fill=black, draw=black] (14,6) rectangle (16, 8);
\filldraw[fill=black, draw=black] (16,2) rectangle (18, 4);
\draw[line width=1pt] (20,0) rectangle (28,8);
\draw[line width=1pt] (22, 0) -- (22,8);
\draw[line width=1pt] (24, 0) -- (24,8);
\draw[line width=1pt] (26, 0) -- (26,8);
\draw[line width=1pt] (20, 2) -- (28,2);
\draw[line width=1pt] (20, 4) -- (28,4);
\draw[line width=1pt] (20, 6) -- (28,6);
\filldraw[fill=black, draw=black] (20,4) rectangle (22, 6);
\filldraw[fill=black, draw=black] (22,2) rectangle (24, 4);
\filldraw[fill=black, draw=black] (24,0) rectangle (26, 2);
\filldraw[fill=black, draw=black] (26,6) rectangle (28, 8);
\draw[line width=1pt] (30,0) rectangle (38,8);
\draw[line width=1pt] (32, 0) -- (32,8);
\draw[line width=1pt] (34, 0) -- (34,8);
\draw[line width=1pt] (36, 0) -- (36,8);
\draw[line width=1pt] (30, 2) -- (38,2);
\draw[line width=1pt] (30, 4) -- (38,4);
\draw[line width=1pt] (30, 6) -- (38,6);
\filldraw[fill=black, draw=black] (30,6) rectangle (32, 8);
\filldraw[fill=black, draw=black] (32,2) rectangle (34, 4);
\filldraw[fill=black, draw=black] (34,4) rectangle (36, 6);
\filldraw[fill=black, draw=black] (36,0) rectangle (38, 2);
\end{tikzpicture}
\caption{Examples of net patterns}
\label{fig:net_patterns}
\end{center}
\end{figure}

 Originally, the net sets were defined \cite{netsets} with the help of \emph{net matrices}, that are $4 \times 4$ matrices having all entries from the set $\{0,1\}$.
Here we give a (shorter) equivalent definition of the net sets, that
 is  appropriate for the framework of the present paper.
We call a $4 \times 4$ pattern a \emph{net pattern} if
each of the four
columns and each of the four rows (each containing four squares) contains exactly one black square, and inside 
each of
the four subsquares of side-length $\frac{1}{2}$ of the unit square that share one vertex with the unit square 
there lies
exactly one of the black squares mentioned above. There exist $16$ such net patterns, four of them are shown in Figure \ref{fig:net_patterns}, the other $12$ can be obtained by flipping or rotating these.
The union of all (closed) white squares in a given net pattern is the corresponding net set of level $0$.

The iterative construction of a sequence of nested sets (\emph{net sets of level} $1,2,\dots$) is analogous to that described in the introduction of Section \ref{sec:introduction} for the Sierpi\'nski carpets, but there is one essential difference: each white square of some 
level can be replaced, by a so-called \emph{net substitution}, by any net set of level $0$ scaled correspondingly, 
such that different white squares can be replaced by different net sets. 
We call a net substitution {\it uniform} if all white squares of a net set of some level are
 substituted by the same scaled net set. Correspondingly, a net set of some level $k\ge 1$ is uniform if at each step of its construction a uniform net substitution (not necessarily the same) was applied.

Thus, by starting with a net pattern and the corresponding net set $E_0$ of level $0$, one obtains, by applying net substitutions,
a decreasing sequence of net sets $E_{0}\supset  E_{1}\supset\dots E_{k-1} \supset E_{k}\supset \dots$.
The fractal $E_{\infty}:=\bigcap_{k\ge 0}E_k$ obtained as the limit set of this construction is called the \emph{limit net set} of the sequence $\{E_k\}_{k\ge 0}$
and can be viewed as the limit set of a Moran construction
\cite{Moran1946, PesinWeiss1996, MauldinWilliams_random_1986} with Hausdorff and box-counting dimension $1+\frac{\log 3}{\log 4}$.

If at each step of the construction we apply a uniform net substitution, not necessarily the same, then
the sets $E_k$ of the above sequence are called \emph{uniform net sets}, 
and $\Elim$ is called
a {\it uniform limit net set}. A \emph{totally uniform net set} is obtained if all substitutions use the same pattern, i.e., we apply the same net substitution at all steps, for all white squares. In this case the obtained limit net set is self-similar: it is  a Sierpi\'nski carpet (as defined in the introduction).

\par When studying connectedness properties of net sets and limit net sets, it is essential to identify two types of patterns: connected and disconnected net patterns. In terms of the graph of the pattern a connected net pattern is a pattern whose graph is connected. Otherwise, the net pattern is called disconnected.
\par Among other, it was proven that connected net patterns (or, originally, connected net
matrices \cite{netsets}) always produce connected limit net sets. For example, in Figure \ref{fig:net_patterns}, the first, second and fourth net pattern are connected, while the third is a disconnected net pattern.

Criteria for different ``degrees'' of connectedness of these fractals have been proven. 
There are four different
possible connectedness ``degrees'' for limit net sets: net-connectedness (a notion introduced in \cite{netsets}), connectedness, disconnectedness and total
disconnectedness \cite{netsets}. It was shown how the connectedness or disconnectedness of the net patterns involved in the iterative construction affects the connectedness ``degree'' of the resulting fractal. 
Necessary and sufficient conditions for the net-connectedness of the fractal were proven, as well as necessary and sufficient conditions for the total disconnectedness of the fractals, necessary and sufficient conditions for a uniform limit net set 
to be connected, but not totally disconnected, or connected (in the Euclidean sense), but not net-connected.

Moreover, an analogon of fractal percolation in the unit cube (see, e.g., \cite{Mandelbrot1983,Falconerbook1990,Falcrandomfr1986}), called \emph{net percolation}, has been introduced, and a sufficient condition for net percolation was proven 
\cite{netsets}. 

The results obtained for limit net sets provide methods for the construction
 of random fractals with a certain type of ``well distributed'' structure (holes) by using net patterns/net matrices, but also for constructing percolating fractal sets and sets that have certain
connectedness properties.

Later on this idea of identifying families of patterns according to their shape was used \cite{totally_disco} in the study of the generalised Sierpi\'nski carpets to which the next section is dedicated.

\section{Generalised Sierpi\'nski carpets}\label{sec:GSC}
  \emph{Generalised Sierpi\'nski carpets} are planar sets in the unit square that were introduced and studied in \cite{connected_general_carpets, totally_disco}. These sets generalise the Sierpi\'nski carpets mentioned in the introduction. They differ in several aspects from a Sierpi\'nski carpet defined as above:  on the one hand, instead of using a single generating pattern, here we use a sequence of patterns in order to construct the generalised carpet,  on the other hand, at any step $k$ of the construction, a $m_k \times m_k$ pattern is used, where $m_k \ge 2$, for all $k \ge 1$, and, moreover, at any two steps $k_1\ne k_2$ we may have distinct patterns, with $m_{k_1} \ne m_{k_2}.$ Thus, generalised Sierpi\'nski carpets are in general not self-similar.
With the notations from Section \ref{sec:patterns} we introduce the following notions.

Let $\{{\cal A}_k\}_{k=1}^{\infty}$ be a sequence of non-empty patterns and $\{m_k\}_{k=1}^{\infty}$ be the corresponding \emph{width-sequence}, i.e., for all $k\ge 1$ we have ${\cal A}_k\subseteq {\cal S}_{m_k}$. We let ${\cal W}_1={\cal A}_{1}$, and call it the \emph{set of white squares of level $1$}. 
 For $n\ge 2$ we define the \emph{set of white squares of level $n$} by 
\[ {\cal W}_n=\bigcup_{W\in {\cal A}_{n}, W_{n-1}\in {\cal W}_{n-1}}\{ P_{W_{n-1}}(W)\}.\]

\begin{figure}[hhhh]\label{W2}
\begin{center}
\begin{tikzpicture}[scale=.18]
\draw[line width=1pt] (0,0) rectangle (20,20);
\draw[line width=0.8pt] (5, 0) -- (5,20);
\draw[line width=0.8pt] (10, 0) -- (10,20);
\draw[line width=0.8pt] (15, 0) -- (15,20);
\draw[line width=0.8pt] (0, 5) -- (20,5);
\draw[line width=0.8pt] (0, 10) -- (20,10);
\draw[line width=0.8pt] (0, 15) -- (20,15);
\draw[line width=0.5pt] (1, 0) -- (1,20);
\draw[line width=0.5pt] (2, 0) -- (2,20);
\draw[line width=0.5pt] (3, 0) -- (3,20);
\draw[line width=0.5pt] (4, 0) -- (4,20);
\draw[line width=0.5pt] (6, 0) -- (6,20);
\draw[line width=0.5pt] (7, 0) -- (7,20);
\draw[line width=0.5pt] (8, 0) -- (8,20);
\draw[line width=0.5pt] (9, 0) -- (9,20);
\draw[line width=0.5pt] (11, 0) -- (11,20);
\draw[line width=0.5pt] (12, 0) -- (12,20);
\draw[line width=0.5pt] (13, 0) -- (13,20);
\draw[line width=0.5pt] (14, 0) -- (14,20);
\draw[line width=0.5pt] (16, 0) -- (16,20);
\draw[line width=0.5pt] (17, 0) -- (17,20);
\draw[line width=0.5pt] (18, 0) -- (18,20);
\draw[line width=0.5pt] (19, 0) -- (19,20);
\draw[line width=0.5pt] (0, 1) -- (20,1);
\draw[line width=0.5pt] (0, 2) -- (20,2);
\draw[line width=0.5pt] (0, 3) -- (20,3);
\draw[line width=0.5pt] (0, 4) -- (20,4);
\draw[line width=0.5pt] (0, 6) -- (20,6);
\draw[line width=0.5pt] (0, 7) -- (20,7);
\draw[line width=0.5pt] (0, 8) -- (20,8);
\draw[line width=0.5pt] (0, 9) -- (20,9);
\draw[line width=0.5pt] (0, 11) -- (20,11);
\draw[line width=0.5pt] (0, 12) -- (20,12);
\draw[line width=0.5pt] (0, 13) -- (20,13);
\draw[line width=0.5pt] (0, 14) -- (20,14);
\draw[line width=0.5pt] (0, 16) -- (20,16);
\draw[line width=0.5pt] (0, 17) -- (20,17);
\draw[line width=0.5pt] (0, 18) -- (20,18);
\draw[line width=0.5pt] (0, 19) -- (20,19);
\filldraw[fill=black, draw=black] (5,0) rectangle (10, 5);
\filldraw[fill=black, draw=black] (10,0) rectangle (15, 5);
\filldraw[fill=black, draw=black] (15,0) rectangle (20, 5);
\filldraw[fill=black, draw=black] (10,5) rectangle (15, 10);
\filldraw[fill=black, draw=black] (0,10) rectangle (5, 15);
\filldraw[fill=black, draw=black] (10,15) rectangle (15, 20);
\filldraw[fill=black, draw=black] (15,15) rectangle (20, 20);
\filldraw[fill=black, draw=black] (1,0) rectangle (2, 1);
\filldraw[fill=black, draw=black] (3,0) rectangle (4, 1);
\filldraw[fill=black, draw=black] (4,0) rectangle (5, 1);
\filldraw[fill=black, draw=black] (4,1) rectangle (5, 2);
\filldraw[fill=black, draw=black] (0,2) rectangle (1, 3);
\filldraw[fill=black, draw=black] (2,2) rectangle (3, 3);
\filldraw[fill=black, draw=black] (2,3) rectangle (3, 4);
\filldraw[fill=black, draw=black] (3,3) rectangle (4, 4);
\filldraw[fill=black, draw=black] (0,4) rectangle (1, 5);
\filldraw[fill=black, draw=black] (3,4) rectangle (4, 5 );
\filldraw[fill=black, draw=black] (4,4) rectangle (5, 5);
\filldraw[fill=black, draw=black] (1,5) rectangle (2, 6);
\filldraw[fill=black, draw=black] (3,5) rectangle (4, 6);
\filldraw[fill=black, draw=black] (4,5) rectangle (5, 6);
\filldraw[fill=black, draw=black] (4,6) rectangle (5, 7);
\filldraw[fill=black, draw=black] (0,7) rectangle (1, 8);
\filldraw[fill=black, draw=black] (2,7) rectangle (3, 8);
\filldraw[fill=black, draw=black] (2,8) rectangle (3, 9);
\filldraw[fill=black, draw=black] (3,8) rectangle (4, 9);
\filldraw[fill=black, draw=black] (0,9) rectangle (1, 10);
\filldraw[fill=black, draw=black] (3,9) rectangle (4, 10);
\filldraw[fill=black, draw=black] (4,9) rectangle (5, 10);
\filldraw[fill=black, draw=black] (6,5) rectangle (7, 6);
\filldraw[fill=black, draw=black] (8,5) rectangle (9, 6);
\filldraw[fill=black, draw=black] (9,5) rectangle (10, 6);
\filldraw[fill=black, draw=black] (9,6) rectangle (10, 7);
\filldraw[fill=black, draw=black] (5,7) rectangle (6, 8);
\filldraw[fill=black, draw=black] (7,7) rectangle (8, 8);
\filldraw[fill=black, draw=black] (7,8) rectangle (8, 9);
\filldraw[fill=black, draw=black] (8,8) rectangle (9, 9);
\filldraw[fill=black, draw=black] (5,9) rectangle (6, 10);
\filldraw[fill=black, draw=black] (8,9) rectangle (9, 10);
\filldraw[fill=black, draw=black] (9,9) rectangle (10, 10);
\filldraw[fill=black, draw=black] (6,10) rectangle (7, 11);
\filldraw[fill=black, draw=black] (8,10) rectangle (9, 11);
\filldraw[fill=black, draw=black] (9,10) rectangle (10, 11);
\filldraw[fill=black, draw=black] (9,11) rectangle (10, 12);
\filldraw[fill=black, draw=black] (7,12) rectangle (6, 13);
\filldraw[fill=black, draw=black] (7,12) rectangle (8, 13);
\filldraw[fill=black, draw=black] (7,13) rectangle (8, 14);
\filldraw[fill=black, draw=black] (8,13) rectangle (9, 14);
\filldraw[fill=black, draw=black] (5,14) rectangle (6, 15);
\filldraw[fill=black, draw=black] (8,14) rectangle (9, 15 );
\filldraw[fill=black, draw=black] (9,14) rectangle (10, 15);
\filldraw[fill=black, draw=black] (11,10) rectangle (12, 11);
\filldraw[fill=black, draw=black] (13,10) rectangle (14, 11);
\filldraw[fill=black, draw=black] (14,10) rectangle (15, 11);
\filldraw[fill=black, draw=black] (14,11) rectangle (15, 12);
\filldraw[fill=black, draw=black] (10,12) rectangle (11, 13);
\filldraw[fill=black, draw=black] (12,12) rectangle (13, 13);
\filldraw[fill=black, draw=black] (12,13) rectangle (13, 14);
\filldraw[fill=black, draw=black] (13,13) rectangle (14, 14);
\filldraw[fill=black, draw=black] (10,14) rectangle (11, 15);
\filldraw[fill=black, draw=black] (13,14) rectangle (14, 15 );
\filldraw[fill=black, draw=black] (14,14) rectangle (15, 15);
\filldraw[fill=black, draw=black] (16,5) rectangle (17, 6);
\filldraw[fill=black, draw=black] (18,5) rectangle (19, 6);
\filldraw[fill=black, draw=black] (19,5) rectangle (20, 6);
\filldraw[fill=black, draw=black] (19,6) rectangle (20, 7);
\filldraw[fill=black, draw=black] (15,7) rectangle (16, 8);
\filldraw[fill=black, draw=black] (17,7) rectangle (18, 8);
\filldraw[fill=black, draw=black] (17,8) rectangle (18, 9);
\filldraw[fill=black, draw=black] (18,8) rectangle (19, 9);
\filldraw[fill=black, draw=black] (15,9) rectangle (16, 10);
\filldraw[fill=black, draw=black] (18,9) rectangle (19, 10 );
\filldraw[fill=black, draw=black] (19,9) rectangle (20, 10);
\filldraw[fill=black, draw=black] (16,10) rectangle (17, 11);
\filldraw[fill=black, draw=black] (18,10) rectangle (19, 11);
\filldraw[fill=black, draw=black] (19,10) rectangle (20, 11);
\filldraw[fill=black, draw=black] (19,11) rectangle (20, 12);
\filldraw[fill=black, draw=black] (15,12) rectangle (16, 13);
\filldraw[fill=black, draw=black] (17,12) rectangle (18, 13);
\filldraw[fill=black, draw=black] (17,13) rectangle (18, 14);
\filldraw[fill=black, draw=black] (18,13) rectangle (19, 14);
\filldraw[fill=black, draw=black] (15,14) rectangle (16, 15);
\filldraw[fill=black, draw=black] (18,14) rectangle (19, 15);
\filldraw[fill=black, draw=black] (19,14) rectangle (20, 15);
\filldraw[fill=black, draw=black] (1,15) rectangle (2, 16);
\filldraw[fill=black, draw=black] (3,15) rectangle (4, 16);
\filldraw[fill=black, draw=black] (4,15) rectangle (5, 16);
\filldraw[fill=black, draw=black] (4,16) rectangle (5, 17);
\filldraw[fill=black, draw=black] (0,17) rectangle (1, 18);
\filldraw[fill=black, draw=black] (2,17) rectangle (3, 18);
\filldraw[fill=black, draw=black] (2,18) rectangle (3, 19);
\filldraw[fill=black, draw=black] (3,18) rectangle (4, 19);
\filldraw[fill=black, draw=black] (0,19) rectangle (1, 20);
\filldraw[fill=black, draw=black] (3,19) rectangle (4, 20);
\filldraw[fill=black, draw=black] (4,19) rectangle (5, 20);
\filldraw[fill=black, draw=black] (6,15) rectangle (7, 16);
\filldraw[fill=black, draw=black] (8,15) rectangle (9, 16);
\filldraw[fill=black, draw=black] (9,15) rectangle (10, 16);
\filldraw[fill=black, draw=black] (9,16) rectangle (10, 17);
\filldraw[fill=black, draw=black] (5,17) rectangle (6, 18);
\filldraw[fill=black, draw=black] (7,17) rectangle (8, 18);
\filldraw[fill=black, draw=black] (7,18) rectangle (8, 19);
\filldraw[fill=black, draw=black] (8,18) rectangle (9, 19);
\filldraw[fill=black, draw=black] (5,19) rectangle (6, 20);
\filldraw[fill=black, draw=black] (8,19) rectangle (9, 20);
\filldraw[fill=black, draw=black] (9,19) rectangle (10, 20);
\filldraw[fill=black, draw=black] (16,15) rectangle (17, 16);
\filldraw[fill=black, draw=black] (18,15) rectangle (19, 16);
\filldraw[fill=black, draw=black] (19,15) rectangle (20, 16);
\filldraw[fill=black, draw=black] (19,16) rectangle (20, 17);
\filldraw[fill=black, draw=black] (15,17) rectangle (16, 18);
\filldraw[fill=black, draw=black] (17,17) rectangle (18, 18);
\filldraw[fill=black, draw=black] (17,18) rectangle (18, 19);
\filldraw[fill=black, draw=black] (18,18) rectangle (19, 19);
\filldraw[fill=black, draw=black] (15,19) rectangle (16, 20);
\filldraw[fill=black, draw=black] (18,19) rectangle (19, 20);
\filldraw[fill=black, draw=black] (19,19) rectangle (20, 20);
\end{tikzpicture}
\caption{The set ${\cal W}_2$, constructed based on the patterns ${\cal A}_1$ and ${\cal A}_2$ shown in Figure \ref{fig:A1A2A3}, that
can also be viewed as a $20$-pattern} \label{fig:W2}
\end{center}
\end{figure}
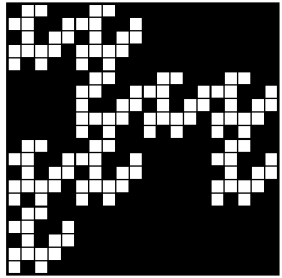

\par By defining three types of graphs associated to the patterns that generate the generalised Sierpi\'nski carpets, necessary and sufficient  conditions for the connectedness
(with respect to the usual
topology of the Euclidean plane) of these objects were proven \cite{connected_general_carpets}.

A different approach, namely identifying certain \emph{families of patterns},
was used \cite{totally_disco} in order to study the structure of the sets obtained at the $n$th iteration in the construction of a generalised carpet, for $n\ge 1$, and it was shown that certain families of patterns provide total disconnectedness of the resulting fractals. Moreover, analogous results hold even in a more general setting \cite{totally_disco}. 
This approach of the carpets provides the possibility to construct disconnected carpets of box-counting dimension less than or even equal to $2$, as it is shown in an example in the more extended arXiv-version  \cite{totally_disco_arxiv} of the published paper \cite{totally_disco}.

The results on connected generalised carpets ~\cite{connected_general_carpets} and  on distances between points on the ``classical'' $s$-dimensional carpet~\cite{carpet_distante_2005}   were extended in more  recent work by Hoffmann~\cite{Hoffmann_GSC_2012}: analogues of the generalised Sierpi\'nski carpets mantioned above, called
\emph{generalised Sierpi\'nski hypercubes}, were defined and studied, and it has been shown that these sets are uniformly
regular, i.e., the geodesic metric is comparable to the Euclidean metric.
We also mention that in previous work,  several authors   
\cite{McMullen_general_carpets_1984, teza_bedford_1984, lalley_gatzouras_1992, peres_selfaffine_carpets_1994, baranski_2007}
studied objects called \emph{general Sierpi\'nski carpets} with respect to dimension and
Hausdorff measure using rectangles instead of squares in the definition. While some of these carpets are self-affine, those defined by Bara\'nski \cite{baranski_2007} are not self-affine.

More recently, geometrical and topological
properties of \emph{fractal squares}, which are self-similar Sierpi\'nski carpets as defined in the introduction, were studied \cite{fractal_squares_teza_roinestad_2010, fractal_squares_hongkong}. We note that the self-similar version of the limit net sets mentioned in Section \ref{sec:limit_net_sets}
and the labyrinth fractals mentioned in Section \ref{sec:laby} are fractal squares. Lau, Luo, and
Rao \cite{fractal_squares_hongkong} study the topological structure of a fractal square by studying the
connected components.
Moreover, recently there is  considerable interest 
to study the Lipschitz equivalence of Cantor sets and of totally disconnected fractal squares, e.g., 
\cite{lipschitz_equiv_Cantor_sets_2012,fractal_squares_2016}. 

In more recent work \cite{mixlaby} dedicated to mixed labyrinth fractals, Steinsky and Cristea gave an other sufficient condition for the total disconnectedness of certain classes of generalised Sierpi\'nski carpets, that occur in relation labyrinth fractals that we present in Section \ref{sec:laby}.

We note the combinatorics-flavoured approach of carpets by identifying special families among the generating patterns that provide certain properties, which was inspired by the results that were obtained for limit net sets, that is different from the approaches of other authors who have studied similar objects. On the other hand, we used graphs a lot as a tool in order to characterise some properties of the patterns or of the prefractals obtained at some (finite) steps of the iterative construction of generalised Sierpi\'nski carpets  \cite{connected_general_carpets} or other carpets \cite{laby4x4, generallaby, mixlaby}, and graphs were also used by other authors, e.g., when dealing with fractal squares \cite{fractal_squares_hongkong}.

\section{Labyrinth fractals}\label{sec:laby}

 An other new family of 
(self-similar) fractals, called \emph{labyrinth fractals},  were introduced and studied during the last decade \cite{laby4x4, generallaby}. These fractal objects are (self-similar) dendrites and a special case of the Sierpi\'nski carpets mentioned at the beginning of Section \ref{sec:introduction}.
First, self similar labyrinth fractals generated by a $4 \times 4$ labyrinth pattern were studied \cite{laby4x4}. Originally, such a generator was called ``labyrinth set''  \cite{laby4x4}, not ``labyrinth pattern'', as used in later work on mixed labyrinth fractals \cite{mixlaby, notearcsmixlaby}. Subsequently, by proving several quite technical lemmas and theorems, the results were extended  \cite{generallaby} to the case of self-similar labyrinth fractals generated by $m \times m$ labyrinth patterns, for any $m \ge 5$. 

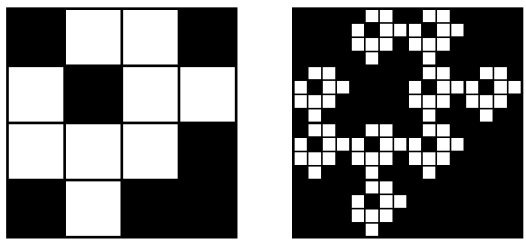
\begin{figure}[hhhh]
\begin{center}
\begin{tikzpicture}[scale=.19]
\draw[line width=1pt] (0,0) rectangle (16,16);
\draw[line width=1pt] (4, 0) -- (4,16);
\draw[line width=1pt] (8, 0) -- (8,16);
\draw[line width=1pt] (12, 0) -- (12,16);
\draw[line width=1pt] (0, 4) -- (16,4);
\draw[line width=1pt] (0, 8) -- (16,8);
\draw[line width=1pt] (0, 12) -- (16,12);
\filldraw[fill=black, draw=black] (0,0) rectangle (4, 4);
\filldraw[fill=black, draw=black] (0,12) rectangle (4, 16);
\filldraw[fill=black, draw=black] (4,8) rectangle (8, 12);
\filldraw[fill=black, draw=black] (8,0) rectangle (12, 4);
\filldraw[fill=black, draw=black] (12,0) rectangle (16, 8);
\filldraw[fill=black, draw=black] (12,12) rectangle (16, 16);
\draw[line width=1pt] (20,0) rectangle (36,16);
\draw[line width=1pt] (24, 0) -- (24,16);
\draw[line width=1pt] (28, 0) -- (28,16);
\draw[line width=1pt] (32, 0) -- (32,16);
\draw[line width=1pt] (20, 4) -- (36,4);
\draw[line width=1pt] (20, 8) -- (36,8);
\draw[line width=1pt] (20, 12) -- (36,12);
\filldraw[fill=black, draw=black] (20,0) rectangle (24, 4);
\filldraw[fill=black, draw=black] (20,12) rectangle (24, 16);
\filldraw[fill=black, draw=black] (24,8) rectangle (28, 12);
\filldraw[fill=black, draw=black] (28,0) rectangle (32, 4);
\filldraw[fill=black, draw=black] (32,0) rectangle (36, 8);
\filldraw[fill=black, draw=black] (32,12) rectangle (36, 16);
\draw[line width=0.6pt] (21, 0) -- (21,16);
\draw[line width=0.6pt] (22, 0) -- (22,16);
\draw[line width=0.6pt] (23, 0) -- (23,16);
\draw[line width=0.6pt] (25, 0) -- (25,16);
\draw[line width=0.6pt] (26, 0) -- (26,16);
\draw[line width=0.6pt] (27, 0) -- (27,16);
\draw[line width=0.6pt] (29, 0) -- (29,16);
\draw[line width=0.6pt] (30, 0) -- (30,16);
\draw[line width=0.6pt] (31, 0) -- (31,16);
\draw[line width=0.6pt] (33, 0) -- (33,16);
\draw[line width=0.6pt] (34, 0) -- (34,16);
\draw[line width=0.6pt] (35, 0) -- (35,16);
\draw[line width=0.6pt] (20, 1) -- (36,1);
\draw[line width=0.6pt] (20, 2) -- (36,2);
\draw[line width=0.6pt] (20, 3) -- (36,3);
\draw[line width=0.6pt] (20, 5) -- (36,5);
\draw[line width=0.6pt] (20, 6) -- (36,6);
\draw[line width=0.6pt] (20, 7) -- (36,7);
\draw[line width=0.6pt] (20, 9) -- (36,9);
\draw[line width=0.6pt] (20, 10) -- (36,10);
\draw[line width=0.6pt] (20, 11) -- (36,11);
\draw[line width=0.6pt] (20, 13) -- (36,13);
\draw[line width=0.6pt] (20, 14) -- (36,14);
\draw[line width=0.6pt] (20, 15) -- (36,15);
\filldraw[fill=black, draw=black] (24,0) rectangle (25, 1);
\filldraw[fill=black, draw=black] (24,3) rectangle (25, 4);
\filldraw[fill=black, draw=black] (25,2) rectangle (26, 3);
\filldraw[fill=black, draw=black] (26,0) rectangle (27, 1);
\filldraw[fill=black, draw=black] (27,0) rectangle (28, 2);
\filldraw[fill=black, draw=black] (27,3) rectangle (28, 4);
\filldraw[fill=black, draw=black] (20,4) rectangle (21, 5);
\filldraw[fill=black, draw=black] (20,7) rectangle (21, 8);
\filldraw[fill=black, draw=black] (21,6) rectangle (22, 7);
\filldraw[fill=black, draw=black] (22,4) rectangle (23, 5);
\filldraw[fill=black, draw=black] (23,4) rectangle (24, 6);
\filldraw[fill=black, draw=black] (23,7) rectangle (24, 8);
\filldraw[fill=black, draw=black] (24,4) rectangle (25, 5);
\filldraw[fill=black, draw=black] (24,7) rectangle (25, 8);
\filldraw[fill=black, draw=black] (25,6) rectangle (26, 7);
\filldraw[fill=black, draw=black] (26,4) rectangle (27, 5);
\filldraw[fill=black, draw=black] (27,4) rectangle (28, 6);
\filldraw[fill=black, draw=black] (27,7) rectangle (28, 8);
\filldraw[fill=black, draw=black] (28,4) rectangle (29, 5);
\filldraw[fill=black, draw=black] (28,7) rectangle (29, 8);
\filldraw[fill=black, draw=black] (29,6) rectangle (30, 7);
\filldraw[fill=black, draw=black] (30,4) rectangle (31, 5);
\filldraw[fill=black, draw=black] (31,4) rectangle (32, 6);
\filldraw[fill=black, draw=black] (31,7) rectangle (32, 8);
\filldraw[fill=black, draw=black] (20,8) rectangle (21, 9);
\filldraw[fill=black, draw=black] (20,11) rectangle (21, 12);
\filldraw[fill=black, draw=black] (21,10) rectangle (22, 11);
\filldraw[fill=black, draw=black] (22,8) rectangle (23, 9);
\filldraw[fill=black, draw=black] (23,8) rectangle (24, 10);
\filldraw[fill=black, draw=black] (23,11) rectangle (24, 12);
\filldraw[fill=black, draw=black] (28,8) rectangle (29, 9);
\filldraw[fill=black, draw=black] (28,11) rectangle (29, 12);
\filldraw[fill=black, draw=black] (29,10) rectangle (30, 11);
\filldraw[fill=black, draw=black] (30,8) rectangle (31, 9);
\filldraw[fill=black, draw=black] (31,8) rectangle (32, 10);
\filldraw[fill=black, draw=black] (31,11) rectangle (32, 12);
\filldraw[fill=black, draw=black] (32,8) rectangle (33, 9);
\filldraw[fill=black, draw=black] (32,11) rectangle (33, 12);
\filldraw[fill=black, draw=black] (33,10) rectangle (34, 11);
\filldraw[fill=black, draw=black] (34,8) rectangle (35, 9);
\filldraw[fill=black, draw=black] (35,8) rectangle (36, 10);
\filldraw[fill=black, draw=black] (35,11) rectangle (36, 12);
\filldraw[fill=black, draw=black] (24,12) rectangle (25, 13);
\filldraw[fill=black, draw=black] (24,15) rectangle (25, 16);
\filldraw[fill=black, draw=black] (25,14) rectangle (26, 15);
\filldraw[fill=black, draw=black] (26,12) rectangle (27, 13);
\filldraw[fill=black, draw=black] (27,12) rectangle (28, 14);
\filldraw[fill=black, draw=black] (27,15) rectangle (28, 16);
\filldraw[fill=black, draw=black] (28,12) rectangle (29, 13);
\filldraw[fill=black, draw=black] (28,15) rectangle (29, 16);
\filldraw[fill=black, draw=black] (29,14) rectangle (30, 15);
\filldraw[fill=black, draw=black] (30,12) rectangle (31, 13);
\filldraw[fill=black, draw=black] (31,12) rectangle (32, 14);
\filldraw[fill=black, draw=black] (31,15) rectangle (32, 16);
\end{tikzpicture}
\caption{A $4 \times 4 $ labyrinth pattern and the corresponding labyrinth set of level 2 (that can also be viewed as a $16 \times 16$ labyrinth pattern)}
\label{fig:laby patterns krabbe}
\end{center}
\end{figure}

To the $m \times m$ pattern that generates a labyrinth fractal we associate the set of ${\cal W}_{1}$ of its (closed) white 
squares and call it the \emph{set of white squares of level 1} or labyrinth set of level 1, and we define $\displaystyle L_1=\cup_{W \in {\cal W}_1}$. By the iterative construction described in Section \ref{sec:introduction} we then obtain the sequence  $\{ {\cal W}_{n}\}_{n \ge 1}$, with ${\cal W}_{n}\subset {{\cal S}_{m^n}}$, and the decreasing sequence of compact sets $\{L_n\}_{n\ge 1}$.

A \emph{top exit} in ${\cal W}_{n}$ is a white
square in the top row of ${\cal W}_{n}$,
 such that there is also a white square in the same column in the bottom row. 
The \emph{bottom exit}, \emph{left exit}, and \emph{right exit} are defined analogously.  



A non-empty $m$-pattern ${\cal A} \subseteq {{\cal S}_m}$, $m \ge 3$ is called a $m\times m$-\emph{labyrinth pattern} (in short, \emph{labyrinth pattern}) if  ${\cal A}$ satisfies satisfies the following three properties:

(1) $\mathcal{G}({\cal W}_{n})$ is a tree;

(2) exactly one top exit in ${\cal W}_{n}$ lies in the top row (of order $n$), exactly one bottom exit lies
in the bottom row, exactly one left exit lies in the left column, and exactly one right exit lies in the right column;

(3) if there is a white square in ${\cal W}_{n}$ at a corner of ${\cal W}_{n}$, then there is no white square in
${\cal W}_{n}$ at the diagonally opposite corner of ${\cal W}_{n}$.

\par
We note that the graph $\mathcal{G}({\cal A})$ introduced in Section \ref{sec:patterns} can also be defined in the case ${\cal A}={\cal W}_{n}$, $n \ge 1$. These graphs play an important role throughout the study of labyrinth sets and labyrinth fractals.

For any labyrinth pattern ${\cal A}$ and any integer $n\ge 1$, the \emph{labyrinth set} (of level $n$) ${\cal W}_n$  
has the above properties  (1), (2), and (3) of a labyrinth pattern \cite{laby4x4}.

The limit set $L_{\infty}$  of the decreasing sequence of compact sets  $\{ L_{n}\}_{n \ge 1}$ is called a \emph{labyrinth fractal}. Every labyrinth fractal has four exits. The top exit of $L_{\infty}$ lies on the top edge of the unit square and is the intersection (point) $\cap_{n=1}^{\infty} T_n$, where $T_n$ is the top exit in ${\cal W}_{n}$, for all $n\ge 1$. The bottom exit, the left and the right exit of $L_{\infty}$ are defined analogously and lie correspondingly on the other edges of the unit square.

We say that an $m\times m$-labyrinth pattern ${\cal A}$  is \emph{horizontally blocked} if the row (of squares) from 
the left to the right exit of ${\cal A}$ contains at least one black square, and it is called \emph{vertically blocked} if the 
column (of squares) from the top to the bottom exit contains at least one black square. 
We remark that in Figure \ref{fig:A1A2A3} the first two labyrinth patterns  are both horizontally and vertically blocked, and the third pattern is neither horizontally, nor vertically blocked, while the $4 \times 4$ pattern shown in Figure \ref{fig:laby patterns krabbe} is both vertically and horizontally blocked. For more examples we refer to \cite{laby4x4,generallaby,mixlaby,notearcsmixlaby}.

Both topological and geometrical properties and aspects of the labyrinth sets and fractals were studied. 
It was proven  that any self-similar labyrinth fractal $L_{\infty}$ is a dendrite, i.e., a locally connected continuum that contains no simple closed curve \cite{laby4x4, generallaby}. 

Subsequently, the arcs in $L_{\infty}$ that connect exits of the fractal were studied, with emphasis on their length. In order to obtain results for the lengths of such arcs, we studied the lengths of paths in the tree $ {\cal G}({\cal W}_n)$ between exits in ${\cal W}_n$. Therefore, the \emph{path matrix} $M$ of the labyrinth set ${\cal W}_1$ was introduced, which is a $6 \times 6$ matrix where each entry represents the number of a certain type of squares in one of the $6$ paths in $ {\cal G}({\cal W}_1)$ between two exits of  ${\cal W}_1$. (For more details on the possible $6$ types of squares in a path in $ {\cal G}({\cal W}_1)$ see, e.g., \cite{laby4x4}.) The path matrix plays an essential role and is a powerful instrument when dealing with lengths of paths in $ {\cal G}({\cal W}_n)$ and with lengths of arcs between exits in $L_{\infty}$. Moreover, this matrix actually is the matrix of a substitution.

In order to prove the obtained results, several known theorems from different areas of mathematics were used: the Perron-Frobenius Theorem, the Hahn-Mazurkiewicz-Sierpi\'nski Theorem (that characterises local connectedness),  the Jordan Curve Theorem, and a labyrinth version of the Steinhaus Chessboard Theorem, proven in \cite{laby4x4}.

It was essential to establish \cite{laby4x4, generallaby} a recursion and prove that the $n$-th power $M^n$ of the path matrix gives information about the lengths of paths in $ {\cal G}({\cal W}_n)$. Moreover, it was shown that the path matrix of a  $m \times m$- labyrinth pattern (or set) is primitive if and only if the pattern (set) is horizontally and vertically blocked. Then, the Perron-Frobenius Theorem for primitive matrices was used, see e.g. \cite[Theorem 1.1, p.3]{Seneta}, in order to obtain the asymptotics of the path lengths in $ {\cal G}({\cal W}_n)$ as $n$ tends to infinity. This subsequently lead to results about the lengths of arcs in the labyrinth fractals generated by both horizontally and vertically blocked labyrinth patterns.
In the case of $m \times m$-patterns with $m \ge 5$, not just the path matrix mentioned above, but also use a second matrix, the reduced path matrix, was used in order to prove the infinite length of arcs in the fractal \cite{generallaby}.

The main results on labyrinth fractals, both in the case when the fractal is generated by a $4 \times 4$ pattern and in the case when when the generating pattern is $m \times m$, with $m\ge 5$, are contained in the following theorem \cite{laby4x4, generallaby}.\\[0.3cm]
{\bf Theorem.} \emph{If $L_{\infty}$ is the labyrinth fractal generated by a horizontally and 
vertically blocked $m\times m$-labyrinth pattern ($m\ge 4$) with path matrix $M$, and $r$ is the spectral radius of $M$, then 
between any two points in $L_{\infty}$ there is a unique arc $\boldsymbol{a}$, the length of $\boldsymbol{a}$ is infinite, and the set of 
all points, at which no tangent to $\boldsymbol{a}$ exists, is dense in $\boldsymbol{a}$.
Moreover, if $\boldsymbol{a}$ is an arc between two distinct points in $L_{\infty}$ then its box-counting dimension is} 
$\dim_B(\boldsymbol{a})=\frac{\log(r)}{\log(m)}$.
\\[0.3cm]
The case when the labyrinth fractal is generated by a $4 \times 4$-pattern that is blocked only in one direction (e.g., only horizontally, but not vertically blocked) is also interesting: we have proven that then there exist both arcs of finite length and arcs of infinite length in the fractal. Moreover, in this case the box-counting dimension of every arc is $1$, while in the case of a labyrinth fractal generated by a both horizontally and vertically blocked pattern the dimension of such arcs is always strictly greater than $1$. For more details we refer to \cite{laby4x4}.

\emph{Mixed labyrinth fractals} were defined and studied later \cite{mixlaby}, as a generalisation of the self-similar labyrinth fractals mentioned above. Here, the construction is analogous to that of the generalised Sierpi\'nski carpets mentioned in Section \ref{sec:GSC}, with the difference that all patterns that occur throughout the construction are labyrinth patterns. In other words, mixed labyrinth fractals are a special case of genersalised Sierpi\'nski carpets. As an example, Figure \ref{fig:W2} shows the labyrinth set of level 2, ${\cal W}_2$, generated by the patterns ${\cal A}_1$ and ${\cal A}_2$ from Figure \ref{fig:A1A2A3}. 

Mixed labyrinth fractals are, like the self-smilar labyrinth fractals mentioned above, dendrites \cite{mixlaby}, but here things get more complicated when studying lengths of paths in the graphs of mixed labyrinth sets of some level, and a lot more complicated when studying lengths of arcs in the fractal. The methods used in the self-similar case, based on the path matrix, can only be appplied up to a certain point in the reasoning. This is due to the fact that in the self similar case the path matrix of a labyrinth set of level $n$  is just $M^n$, where $M$ is the path matrix of the generating pattern, while in the mixed case it is $M_1 \cdot M_2\cdot M_n$, where $M_k$ is the path matrix of the pattern ${\cal A}_k,$ for $k=1,\dots,n$. 
Since there are no restrictions regarding the labyrinth patterns that occur in the generating sequence $\{{\cal A_k}\}_{k\ge 1}$, the methods used in the self-similar case in order to establish results on the asymptotical behavior of the path matrix associated to the labyrinth set of level $n$, for $n \to \infty,$ cannot be applied here anymore, as soon as the sequence contains more than one pattern (unless the sequence is periodic, which is in general not the case). Moreover, there are also other properties that get lost when we give up self-similarity, e.g., in the case of mixed labyrinth fractals it is possible that an exit of the fractal lies, for some $n\ge 1$, in more than one white square of level $n$ of ${\cal W}_n$, while in the self-similar case each exit lies in a unique such white square of level $n$.

 In a recently published paper \cite{notearcsmixlaby} in was shown, that in the case of a mixed labyrinth fractal the theorem stated above does not hold. More precisely, one can prove the following two results.
\\
{\bf Theorem.} \emph{
There exist  sequences $\{{\cal A}_k\}_{k=1}^{\infty}$ of (both horizontally and vertically) blocked labyrinth patterns, such that the limit set $L_{\infty}$ has the property that
 for any two points in $L_{\infty}$ the length of the arc $a\subset L_{\infty}$ that connects them is finite. For almost all points $x_0 \in a$ (with respect to the length) there exists the tangent at $x_0$ to the arc $a$.}
\\[0.1cm]
{\bf Proposition.} \emph{There exist  sequences $\{{\cal A}_k\}_{k=1}^{\infty}$ of (both horizontally and vertically) blocked labyrinth patterns, such that the limit set $L_{\infty}$ has the property that
 for any two points in $L_{\infty}$ the length of the arc $a\subset L_{\infty}$ that connects them is infinite.}

\begin{figure}[h!]
\begin{center}
\begin{tikzpicture}[scale=.15]
\draw[line width=1.5pt] (0,0) rectangle (22,22);
\draw[line width=1pt] (2, 0) -- (2,22);
\draw[line width=1pt] (4, 0) -- (4,22);
\draw[line width=1pt] (6, 0) -- (6,22);
\draw[line width=1pt] (8, 0) -- (8,22);
\draw[line width=1pt] (10, 0) -- (10,22);
\draw[line width=1pt] (12, 0) -- (12,22);
\draw[line width=1pt] (14, 0) -- (14,22);
\draw[line width=1pt] (16, 0) -- (16,22);
\draw[line width=1pt] (18, 0) -- (18,22);
\draw[line width=1pt] (20, 0) -- (20,22);
\draw[line width=1pt] (0, 2) -- (22,2);
\draw[line width=1pt] (0, 4) -- (22,4);
\draw[line width=1pt] (0, 6) -- (22,6);
\draw[line width=1pt] (0, 8) -- (22,8);
\draw[line width=1pt] (0, 10) -- (22,10);
\draw[line width=1pt] (0, 12) -- (22,12);
\draw[line width=1pt] (0, 14) -- (22,14);
\draw[line width=1pt] (0, 16) -- (22,16);
\draw[line width=1pt] (0, 18) -- (22,18);
\draw[line width=1pt] (0, 20) -- (22,20);
\filldraw[fill=black, draw=black] (0,0) rectangle (8, 10);
\filldraw[fill=black, draw=black] (8,8) rectangle (10, 10);
\filldraw[fill=black, draw=black] (8,0) rectangle (10, 2);
\filldraw[fill=black, draw=black] (4,10) rectangle (6, 12);
\filldraw[fill=black, draw=black] (0,14) rectangle (10, 22);
\filldraw[fill=black, draw=black] (0,12) rectangle (2, 17);
\filldraw[fill=black, draw=black] (8,12) rectangle (10, 14);
\filldraw[fill=black, draw=black] (10,16) rectangle (12, 18);
\filldraw[fill=black, draw=black] (12,20) rectangle (14, 22);
\filldraw[fill=black, draw=black] (12,12) rectangle (14, 14);
\filldraw[fill=black, draw=black] (14,12) rectangle (22, 22);
\filldraw[fill=black, draw=black] (12,8) rectangle (14, 10 );

\filldraw[fill=black, draw=black] (16,10) rectangle (18, 12);
\filldraw[fill=black, draw=black] (12,0) rectangle (22, 8);
\filldraw[fill=black, draw=black] (20,8) rectangle (22, 10);
\filldraw[fill=black, draw=black] (10,4) rectangle (12, 6);
\end{tikzpicture}
\caption{An example: a special cross pattern with width $11$}
\label{fig:specialcross}
\end{center}
\end{figure}
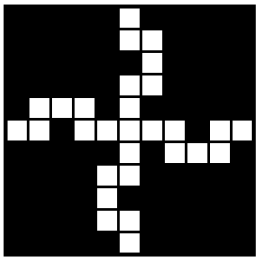

These results were proven by using a special family of labyrinth patterns, which we called ``special cross patterns''. An example of such a pattern is shown in Figure 
\ref{fig:specialcross}.

Here, the idea was to approximate the arcs that connect exits in the labyrinth fractal by special curves, which are related to the patterns used in the construction. For more details we refer to the paper \cite{notearcsmixlaby}. Analogously to the self-similar case, by chosing suitable labyrinth patterns that are blocked in only one direction (e.g., horizontally but not vertically blocked), one can construct also in the mixed case labyrinth fractals where both arcs of finite length and arcs of infinite length exist between their points.

The following conjecture was formulated \cite{notearcsmixlaby}:\\
{\bf Conjecture.}  \emph{A sequence of both horizontally and vertically blocked labyrinth patterns with the property that the sequence of widths $\{m_k\}_{k\ge 1}$ is bounded, generates a mixed labyrinth fractal with the property that for any $x,y \in L_{\infty}$ the length of the arc in the fractal that connects $x$ and $y$ is infinite.}

\begin{figure}
\begin{center}
 \begin{tikzpicture}[scale=.25]

\draw[line width=1pt] (0,0) rectangle (14,14);
\draw[line width=1pt] (2, 0) -- (2,14);
\draw[line width=1pt] (4, 0) -- (4,14);
\draw[line width=1pt] (6, 0) -- (6,14);
\draw[line width=1pt] (8, 0) -- (8,14);
\draw[line width=1pt] (10, 0) -- (10,14);
\draw[line width=1pt] (12, 0) -- (12,14);
\draw[line width=1pt] (0, 2) -- (14,2);
\draw[line width=1pt] (0, 4) -- (14,4);
\draw[line width=1pt] (0, 6) -- (14,6);
\draw[line width=1pt] (0, 8) -- (14,8);
\draw[line width=1pt] (0, 10) -- (14,10);
\draw[line width=1pt] (0, 12) -- (14,12);
\filldraw[fill=black, draw=black] (12,0) rectangle (14, 2);
\filldraw[fill=black, draw=black] (2,2) rectangle (6, 4);
\filldraw[fill=black, draw=black] (8,2) rectangle (10, 6);
\filldraw[fill=black, draw=black] (0,4) rectangle (4, 6);
\filldraw[fill=black, draw=black] (0,8) rectangle (4, 14);
\filldraw[fill=black, draw=black] (10,4) rectangle (14, 6);
\filldraw[fill=black, draw=black] (4,10) rectangle (6, 14);
\filldraw[fill=black, draw=black] (6,6) rectangle (8, 8);
\filldraw[fill=black, draw=black] (10,8) rectangle (14, 10);
\filldraw[fill=black, draw=black] (8,10) rectangle (14, 14);

\hspace{0.4cm}
\draw[line width=1pt] (15,0) rectangle (27,12);
\draw[line width=1pt] (17, 0) -- (17,12);
\draw[line width=1pt] (19, 0) -- (19,12);
\draw[line width=1pt] (21, 0) -- (21,12);
\draw[line width=1pt] (23, 0) -- (23,12);
\draw[line width=1pt] (25, 0) -- (25,12);
\draw[line width=1pt] (15, 2) -- (27,2);
\draw[line width=1pt] (15, 4) -- (27,4);
\draw[line width=1pt] (15, 6) -- (27,6);
\draw[line width=1pt] (15, 8) -- (27,8);
\draw[line width=1pt] (15, 10) -- (27,10);
 \filldraw[fill=black, draw=black] (15,0) rectangle (17, 4);
 \filldraw[fill=black, draw=black] (19,0) rectangle (27, 2);
\filldraw[fill=black, draw=black] (23,2) rectangle (27, 4);
 \filldraw[fill=black, draw=black] (19,4) rectangle (21, 6);
 \filldraw[fill=black, draw=black] (15,6) rectangle (17, 8);
 \filldraw[fill=black, draw=black] (23,6) rectangle (25, 10);
 \filldraw[fill=black, draw=black] (15,8) rectangle (19, 10);
 \filldraw[fill=black, draw=black] (21,10) rectangle (27, 12);
 \end{tikzpicture}
\end{center}
\caption{Examples: two wild labyrinth patterns, both vertically and horizontally blocked}
\label{fig:wild_patterns}
\end{figure}
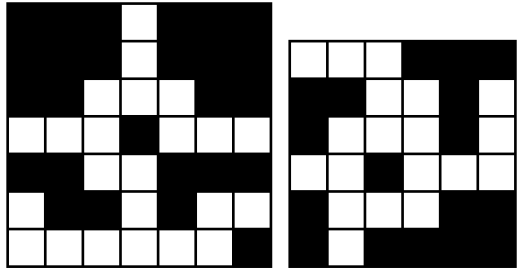

Finally, let us mention that  in \cite{mixlaby} it is also shown how, by relaxing the conditions imposed on the labyrinth patterns in order to construct \emph{wild labyrinth patterns} and, correspondingly, (self-similar or mixed) \emph{wild labyrinth fractals},  several properties of the labyrinth sets and fractals change, and in general the path matrices can not 
be used anymore in order to provide reliable information and results about the paths in the graphs of the labyrinth sets ${\cal G}({\cal W}_n)$, for $n \ge 2$ or the lengths of arcs in the fractal. Figure \ref{fig:wild_patterns} shows two wild labyrinth patterns. The first one has more than one horizontal exit pair, and the graph of the second one is only connected, but not a tree. Moreover, the connected net patterns mentioned in Section \ref{sec:limit_net_sets} are in particular $4 \times 4$ wild labyrinth patterns. 

\section{Conclusions}\label{sec:conclusions}
Although the fractals presented in this paper, limit net sets, generalised Sierpi\'nski carpets and labyrinth fractals, can also be studied by using other approaches than the one used here, e.g., under the framework of IFS, graph-directed constructions,  GDMS, random fractals  or $V$-variable fractals, here we chose this unifying and rather combinatorial approach based on planar patterns. The motivation of this fact is that, in our opinion, this is a way to bring them closer to other sciences, to specialists from other fields, to a wider audience in general, since the notion of ``pattern'' is very intuitive, wide-spread (even in every day life) and somehow basic for the understanding. Moreover, as it follows from the cited papers, patterns are sufficient in order to construct and study these new, special families of fractals with remarkable properties.

The idea is that in this approach mainly by identifying families of patterns or a few properties of the patterns that are easy to check, one can generate fractals of prefractals with desired topological or geometrical properties (like types or degrees of connectedness, or lengths of arcs between points in the fractal), or, in addition, with desired fractal dimension.

It is important to add here that physicists use fractals like those mentioned above as models in different areas, e.g, for the study of materials or of diffusion in porous matter \cite{the_pore_structure_2001,tarafdar_modelling_porous_structures_2001,laby_fizicieni1}, planar nanostructures \cite{laby_fizicieni2}, or even for the construction of new, more performant devices (e.g. radar antennas) \cite{PotapovZhang_oct2016}. 

In this context it is worth to remark that, while the mathematicians focus mainly on the fractal, i.e., on the objects obtained as the limit of the iterative construction, the physicists are usually more interested in some prefractal obtained after a high enough number of iterations. In other words, while the mathematicians are mainly interested in what happens in the infinite, in the limit, the physicsts are interested in what is obtained at a finite step that approximates the limit well enough, or where the scale is fine enough, but not infinitely fine. Let us give an example. From the point of view of research in physics, the properties (1) and (2) of labyrinth patterns are essential, and they are sufficient, since the property (3) of labyrinth patterns only plays a role in the limit, when dealing with the resulting labyrinth fractal: the fact that the resulting labyrinth fractal is a dendrite played an essential role when proving some of the other results, e.g., those on the length of arcs in the fractal.

Finally, let us mention that due to the interesting properties of limit net sets, generalised Sierpi\'nski carpets and labyrinth fractals, at the moment this research is continued on new families of fractals, that are further generalisations of the objects mentioned here: there is a lot of magic and still a lot to discover in this field.

\end{document}